\documentclass[12pt]{article}
\usepackage{amsxtra,amssymb,amsthm,amsmath,latexsym}

\theoremstyle{plain}
\def\R{{\mathbb R}}

\def\oH{{\overset{\circ}{H}}}
\def\oH1{{\overset{\circ}{H}\kern-.02in{}^1}}

\def\bee{\begin{equation*}}
\def\eee{\end{equation*}}
\def\be{\begin{equation}}
\def\ee{\end{equation}}

\begin{document}

%\begin{titlepage}
\title{Global existence and estimates of the solutions to nonlinear integral equations
}

\author{Alexander G. Ramm\\
 Department  of Mathematics, Kansas State University, \\
 Manhattan, KS 66506, USA\\
ramm@math.ksu.edu\\
%,\\ %fax 785-532-0546, tel. 785-532-0580}
http://www.math.ksu.edu/\,$\widetilde{\ }$\,ramm
}

\date{}
\maketitle\thispagestyle{empty}

%%%%%%%%%%%%%%%%%%%%%%%%%%%%%%%%%%%%%%%%%%%%%%%%%
\begin{abstract}

It is proved that a class of nonlinear integral equations of the Volterra-Hammerstein type has a global
solution, that is, solutions defined for all $t\ge 0$,  and estimates of these solutions as $t\to \infty$ are obtained.
The argument uses a nonlinear differential inequality which was proved by the author and has broad
applications.
\end{abstract}
%%%%%%%%%%%%%%%%%%%%%%%%%%%%%%%%%%%%%%%%%%%%%%%%%
%\end{titlepage}

\section{Introduction}\label{S:1}
Consider the equation:
\be\label{e1}
u(t)=\int_0^t e^{-a(t-s)} h(u(s))ds +f(t):=T(u), \quad t\ge 0; \quad a=const>0.
\ee
that is, Volterra-Hammerstein equation. There is a large literature on nonlinear integral equations, \cite{Z}, \cite{D}.
The usual methods to study such equations include fixed-point theorems such as contraction mapping principle and degree theory, (Schauder and Leray-Schauder theorems).
The goal of this paper is to give a new approach to a study of equation  \eqref{e1}. We give sufficient conditions for the global existence of solutions to  \eqref{e1} and their estimates as $t\to \infty$.

Denote $f':=\frac {df}{dt}$. By $c>0$ various constants will be denoted.

Let us formulate our assumptions:
\be\label{e2}
|h(u)|\le c|u|^b, \quad |h'(u)|\le c|u|^{b-1}, \quad b\ge 2,
\ee
\be\label{e3}
|f(t)|+a|f'(t)|\le ce^{-a_1t}, \qquad a_1=const>0.
\ee
By $c>0$ various constants are denoted.

Our approach is based on the author's results on the nonlinear differential inequality
formulated in Theorem 1 (see \cite{R570}--\cite{R657}). These results have been used by the author in a study
of stability of solutions to abstract nonlinear evolution problems (\cite{R657}).

Denote $\R_+=[0,\infty)$.

{\bf Theorem 1.}  {\em Let $g\ge 0$ solve the inequality
\be\label{e4}
g'(t)\le -ag(t) +\alpha(t,g) +\beta(t), \quad t\ge 0, \quad a=const>0,
\ee
where $\alpha(t,g)\ge 0$ and $\beta(t)\ge 0$ are continuous functions of $t$, $t\in \R_+$ and
$\alpha(t,g)$ is locally Lipschitz with respect to $g$. If there exists a function $\mu(t)>0$,
defined on $ \R_+$, $\mu\in C^1( \R_+)$, such that
\be\label{e5}
\alpha(t, \frac 1 {\mu(t)})+\beta(t)\le  \frac 1 {\mu(t)}\Big(a- \frac{\mu'(t)}{\mu(t)}\Big), \quad \forall t\ge 0,
\ee
and
\be\label{e6}
g(0)\mu(0)\le 1,
\ee
then $g$ exists on $\R_+$ and
\be\label{e7}
0\le g(t)\le \frac 1 {\mu(t)}, \quad \forall t\ge 0.
\ee
 }
\vspace{3mm}
A proof of Theorem 1 can be found in \cite{R657}. Its idea is described in Section 2.

The result of this paper is formulated in Theorem 2.

{\bf Theorem 2.}  {\em Assume  that \eqref{e2} and  \eqref{e3} hold, $a\ge 2$, $b\ge 2$, $c\in (0, 0.75)$, $p\in (0, \min(0.75 a, a_1))$, $R=(b-1)^{1/b}$. Then any solution to \eqref{e1} exists on $\R_+$ and satisfies the estimate
\be\label{e8}
|u(t)|\le R^{-1}e^{-pt}, \qquad \forall t\ge 0, \qquad p\in (0, min(0.25 a_1, a)).
\ee
}

In Section 2   Theorem 2 is proved.

\section{ Proof of Theorem 2}\label{S:2}

Let us reduce equation \eqref{e1} to the form suitable for an application of Theorem 1.
Differentiate \eqref{e1} and get
\be\label{e9}
u'=f'-a\int_0^t e^{-a(t-s)} h(u(s))ds +h(u(t)).
\ee
Let $g(t):=|u(t)|$ and take into account that $|F(t)|\le ce^{a_1t}$.

From \eqref{e1}  one gets $\int_0^t e^{-a(t-s)} h(u(s))ds=u-f$. This and equation \eqref{e9} imply
$u'=f'-a(u-f)+h(u(t))$. Therefore, one gets
\be\label{e10}
u'=-au +h(u) +F, \quad F:=f'+af
\ee
Multiply \eqref{e10} by $\overline{u}$, where  $\overline{u}$ stands for complex conjugate of $u$, and get
\be\label{e11}
u'\overline{u}=-ag^2+h(u)\overline{u} +F\overline{u}.
\ee
One has
\be\label{e12}
u'\overline{u}+u(\overline{u})'=\frac {d g^2}{dt}=2gg'.
\ee
We define the derivative as $g'=\lim_{h\to +0} \frac {g(t+h)-g(t)}{h}$. With this
definition,  $g(t)$ is differentiable at every point if $u(t)$ is continuously differentiable
for all $t\ge 0$.
Any solution $u(t)$ to  \eqref{e1} is continuously differentiable under our assumptions.
Take complex conjugate of \eqref{e11}, add  the resulting equation to \eqref{e11}
and take into account \eqref{e12}. This yields
\be\label{e13}
2gg'=-2ag^2 +2 Re (h(u)\overline{u})+ 2 Re (F\overline{u}).
\ee
Since $g\ge 0$,  one derives from   \eqref{e13}, using assumptions \eqref{e2} and \eqref{e3}, that
\be\label{e14}
g'(t)\le -ag(t)+ cg^b+ce^{-a_1t}.
\ee
Let
\be\label{e15}
\mu(t)=Re^{pt}, \quad R=const>0, \quad p\in (0, \min(0.25 a, a_1)).
\ee
Condition \eqref{e5} can be written as
\be\label{e16}
\frac c{R^{b}e^{bpt}} +c e^{-a_1t}\le \frac 1 {Re^{pt}} (a-p), \quad t\in \R_+.
\ee
This inequality holds if
\be\label{e17}
\frac c{R^{b-1}e^{(b-1)pt}} +cR e^{-(a_1-p)t}\le \frac {3a}{4},  \quad t\in \R_+.
\ee
Inequality \eqref{e17} holds if
\be\label{e18}
\frac 1 {R^{b-1}} +R\le \frac {3a}{4c}.
\ee
The minimum of the left side of \eqref{e18} is attained at $R=(b-1)^{1/b}$ and is equal to $\frac b{(b-1)^{(b-1)/b}}$.
Thus,  \eqref{e18} holds if
 \be\label{e19}
 \frac b{(b-1)^{(b-1)/b}}\le \frac {3a}{4c}.
 \ee
For example, assume that
 $$a\ge 2, \qquad c\le 0.75.$$ 
 Then  \eqref{e19} holds if $b\le 2(b-1)^{(b-1)/b}$, that is, if
\be\label{e20}
b^b\le 2^b (b-1)^{b-1}.
\ee
Inequality \eqref{e20} holds if $b\ge 2$. Thus, by Theorem1, any solution $u(t)$ of \eqref{e1} exists globally and
\be\label{e21}
|u(t)|\le \frac {e^{-pt}}{R},
\ee
provided that
\be\label{e22}
|u(0)|R\le 1, \quad R=(b-1)^{1/b}, \quad a\ge 2, \quad b\ge 2, \quad c=0.75, \quad p\in (0, \min(0.25 a, a_1)).
\ee
Inequality $|u(0)|R\le 1$ holds if $f(0)R\le 1$. By assumption \eqref{e3} this inequality holds if $c\le \frac 1 R$.
Theorem 2 is proved. \hfill$\Box$

 Let us prove existence of a solution to \eqref{e1}  using the contraction mapping principle and Theorem 2.

By estimate \eqref{e21} one has $|u(t)|\le \frac 1 R$ for all $t\ge 0$. Therefore, using assumptions \eqref{e2} and \eqref{e3},
one gets
\be\label{e23}
|Tu|\le c+\frac c{aR^b}\le \frac 1 R,
\ee
provided that $cR\le \frac 1{1+\frac 1 {aR^b}}$. For $R=(b-1)^{1/b}$ this inequality holds if $c$ is sufficiently small.
If \eqref{e23} holds, then $T$ maps the ball $B_R:=\{u: ||u||\le \frac 1 R\}$ into itself. Here $||u||=max_{t\ge 0}|u(t)|$.

On the ball $B_R$ the operator $T$ is a contraction:
\be\label{e24}
||Tu-Tv||\le ||\int_0^t e^{-a(t-s)}c|\eta^{b-1}|ds||||u-v||\le \frac {c}{R^{b-1} a} ||u-v||,
\ee
where the assumption \eqref{e2} was used, and $\eta$ is the "middle" element between $u$ and $v$, $||\eta||\le \frac 1 R$.
The integral
in \eqref{e24} is estimated as follows:
\be\label{e25}
||\int_0^t e^{-a(t-s)}c|\eta^{b-1}|ds||\le \frac {c}{R^{b-1}} max_{t\ge 0}\int_0^t  e^{-a(t-s)}ds\le  \frac {c}{R^{b-1} a}
\ee
If
\be\label{e26}
 \frac {c}{R^{b-1} a}<1,
 \ee
 then $T$ is a contraction on $B_R$. Condition \eqref{e26} holds if $c$ is sufficiently small.
Thus, if condition \eqref{e26} and the assumptions of Theorem 2 hold, then, by the contraction mapping principle,
there exists a unique solution to \eqref{e1} in the ball $B_R$. \hfill$\Box$

For convenience of the reader we sketch the idea of the proof of Theorem 1 following \cite{R570}---\cite{R657}.

Inequality \eqref{e5} can be written for the function $w=\frac 1 {\mu}$ as follows:
\be\label{e27}
-a w +\alpha(t,w) +\beta (t)\le w'.
\ee
From \eqref{e4} and  \eqref{e27}  by a comparison lemma for ordinary differential equations it follows that
\be\label{e28}
0\le g(t)\le \frac 1 {\mu(t)},
\ee
provided that $g(0)\le w(0)=\frac 1 {\mu(0)}$. The last inequality is the assumption \eqref{e6}.
Since $\mu(t)>0$ and is assumed to be defined for all $t\ge 0$, the function $w=\frac 1 {\mu}$
is defined for all $t\ge 0$. Since $0\le g(t)\le \frac 1 {\mu(t)}$, and $g(t):=|u(t)|$, the function $u$
is defined for all $t\ge 0$.

If $\lim_{t\to \infty} \mu(t)=0$, then $\lim_{t\to \infty}|u(t)|=0$ by estimate \eqref{e28}.

\end{document}